\documentclass[12pt]{amsart}
\usepackage{amsfonts, amsmath, amsthm}
\usepackage{epsfig}
\usepackage{psfrag}

\newtheorem{pro}{Proposition}[section]
\newtheorem{thm}[pro]{Theorem}
\newtheorem{lem}[pro]{Lemma}

\newtheorem{cor}[pro]{Corollary}

\theoremstyle{definition}
\newtheorem{dfn}[pro]{Definition}
\newtheorem{ex}[pro]{Example}

\theoremstyle{remark}

\title{Normalizing Heegaard-Scharlemann-Thompson Splittings} 
\date{\today}
\address{Mathematics Department, University of Illinois at Chicago}
\email{bachman@math.uic.edu}
\author{David Bachman}
\begin{document}
\begin{abstract}
We define a Heegaard-Scharlemann-Thompson (HST) splitting of a 3-manifold $M$ to be a sequence of pairwise-disjoint, embedded surfaces, $\{F_i\}$, such that for each odd value of $i$, $F_i$ is a Heegaard splitting of the submanifold of $M$ cobounded by $F_{i-1}$ and $F_{i+1}$. Our main result is the following: Suppose $M$ ($\neq B^3$ or $S^3$) is an irreducible submanifold of a triangulated 3-manifold, bounded by a normal or almost normal surface, and containing at most one maximal normal 2-sphere. If $\{F_i\}$ is a strongly irreducible HST splitting of $M$ then we may isotope it so that for each even value of $i$ the surface $F_i$ is normal and for each odd value of $i$ the surface $F_i$ is almost normal. 

We then show how various theorems of Rubinstein, Thompson, Stocking and Schleimer follow from this result. We also show how our results imply the following: (1) a manifold that contains a non-separating surface contains an almost normal one, and (2) if a manifold contains a normal Heegaard surface then it contains two almost normal ones that are topologically parallel to it. 
\end{abstract}
\maketitle


\section{Introduction}
Since Rubinstein's announcement of an algorithm to recognize the 3-sphere in 1993 \cite{rubinstein:93}, almost normal surfaces, the main technical innovation of that work, have received wide attention. In 1994 Thompson \cite{thompson:94} reworked Rubinstein's proof of the existence of an almost normal 2-sphere in $S^3$ using Gabai's language of {\it thin position} \cite{gabai:87}. Also in 1993 Rubinstein announced that any strongly irreducible Heegaard spliting could be made almost normal. This result, combined with his later work with William Jaco \cite{jr:98}, implies that any non-Haken 3-manifold admits at most finitely many Heegaard splittings of a given genus. In 1996 Stocking \cite{stocking:96} provided a new proof of the fact that strongly irreducible Heegaard splittings can be made almost normal, by following Thompson's use of thin position type arguments. 

More recently the author \cite{bachman:98} and, independently, Rieck and Sedgwick \cite{sr:00} used almost normal surfaces to prove that for any orientable 3-manifold with torus boundary, $X$, in all but finitely many Dehn fillings of $X$ the core of the attached solid torus can be isotoped onto every strongly irreducible Heegaard splitting. Saul Schleimer \cite{schleimer:01} has used almost normal surfaces to show that the set of distances of the Heegaard splittings of a given 3-manifold is finite. 

The main result of this paper is that one may isotope a {\it Heegaard-Scharlemann-Thompson} (HST) splitting into an appropriate normal form. Roughly speaking, an HST splitting of a 3-manifold is a sequence of surfaces such that each consecutive pair cobounds a compression body, and if $W$ is such compression body then the subindex of $\partial _+W$ is odd (a more precise definition is given in Section \ref{s:HST}). For each odd $i$ the surface $F_i$ is referred to as a {\it thick level}, and for each even $i$, $F_i$ is a {\it thin level}. Just as for classical Heegaard splittings, the concept of {\it strong irreducibility} is a useful notion of non-triviality for an HST splitting. The precise statement of our main result is the following:

\medskip
\noindent {\bf Theorem \ref{t:main}} {\it Suppose $M$ ($\neq B^3$ or $S^3$) is an irreducible submanifold of a triangulated 3-manifold which is bounded by a normal or almost normal surface, and contains at most one maximal normal 2-sphere. If $\{F_i\}$ is a strongly irreducible HST splitting of $M$ then each thin level is isotopic to a normal surface, and each thick level is isotopic to an almost normal surface.}
\medskip 

Consequences of Theorem \ref{t:main} include many of aforementioned results on the existence of almost normal surfaces. We also show that several new results follow that may be of independent interest. 

\medskip
\noindent {\bf Theorem \ref{nonsep}} {\it If a triangulated 3-manifold contains a non-separating surface then it contains a non-separating almost normal surface.}

\medskip
\noindent {\bf Theorem \ref{compressible}} {\it Suppose $M$ is a submanifold of a triangulated irreducible 3-manifold which is bounded by a connected normal surface, and contains at most one maximal normal 2-sphere. If $\partial M$ is compressible in $M$ then there is an almost normal surface in $M$ which is topologically parallel to $\partial M$, and normalizes monotonically to $\partial M$.}

\medskip

A result similar to Theorem \ref{compressible} has been independently obtained by Culler, Dunfield, Jaco and Shalen in their program to solve the weak Lopez conjecture. An immediate corollary is the following:

\medskip
\noindent {\bf Corollary \ref{normalheegaard}} {\it If $H$ is a Heegaard surface for a 3-manifold which is isotopic to a normal surface in the complement of a maximal normal 2-sphere then $H$ is isotopic to two, topologically parallel, almost normal surfaces.}
\medskip

We conclude the paper by proving a result similar to Theorem \ref{t:main}, which gives an exact characterization of handlebodies in terms of the existence of certain normal and almost normal surfaces. We hope to use this result in subsequent work to give a new algorithm to determine when an almost normal surface is a Heegaard surface.

The techniques used to prove Theorem \ref{t:main} may be of independent interest. While in this paper we make crucial use of ideas from such earlier papers as \cite{thompson:94}, our arguments differ from those of previous authors in several key respects. We outline these differences here. 

\begin{itemize}
\item For the existence of almost normal strongly irreducible Heegaard splittings both Rubinstein and Stocking produce a sequence of surfaces which eventually converge to the desired one. For example, in Rubinstein's  ``iterated" sweepout approach he shows that almost normal surfaces occur at the maximum of a sweepout which is in some way ``minimal". But the first sweepout which he chooses may not have the desired surface as a maximum. So a sequence of sweepouts of smaller and smaller submanifolds is chosen in which the maximum is eventually the desired surface. In the approach here we show that if the right concept of ``minimal" is chosen for the initial sweepout then all {\it local} maxima are almost normal, and the desired surface must be one of these. 

\item Several results are used here to streamline the argument which were not available to Rubinstein or Stocking. For example, we appeal to Scharlemann's {\it Local detection of strongly irreducible Heegaard splittings} \cite{scharlemann:97} to classify the intersection of certain surfaces with the interior of a 3-simplex.

\item Rather than using sweepouts (as does Rubinstein) or thin position (as do Thompson and Stocking) we use the general framework of {\it relative} HST Splittings (see Section \ref{s:relHST}). We find that this general framework provides a much cleaner approach and further streamlines some of the arguments. 

\item We deal with 2-sphere components that may be ``pinched off" in the formation of an almost normal surface in an entirely different way than previous authors. 
\end{itemize}


The author thanks Saul Schleimer for helpful conversations, especially regarding Lemma \ref{l:untang}, during the preparation of this paper. The author also thanks Peter Shalen.

\section{Background material.}
In this section, we give some of the standard definitions that will be used throughout the paper. The expert in 3-manifold theory may still want to skim this, as one or two new terms are introduced. 

A 2-sphere in a 3-manifold which does not bound a 3-ball on either side is called {\it essential}. If a manifold does not contain an essential 2-sphere, then it is referred to as {\it irreducible}.

A loop on a surface is called {\it essential} if it does not bound a disk in the surface. Given a surface, $F$, in a 3-manifold, $M$, a {\it compressing disk} for $F$ is a disk, $D \subset M$, such that $F \cap D=\partial D$, and such that $\partial D$ is essential on $F$. If we let $D \times I$ denote a thickening of $D$ in $M$, then to {\it compress $F$ along $D$} is to remove $(\partial D) \times I$ from $F$, and replace it with $D \times \partial I$.  

A {\it compression body} is a 3-manifold which can be obtained by starting with some surface, $F$ (not necessarily connected), forming the product, $F \times I$, attaching some number of 2-handles to $F \times \{1\}$, and capping off all remaining 2-sphere boundary components with 3-balls. The boundary component, $F \times \{0\}$, is often referred to as $\partial _+$. The other boundary component is referred to as $\partial _-$. If $\partial _-=\emptyset$, then we say the compression body is a {\it handlebody}. A compression body is {\it non-trivial} if it is not a product. 

A surface, $F$, in a 3-manifold, $M$, is a {\it Heegaard surface for M} if $F$ separates $M$ into two compression bodies, $W$, and $W'$, such that $F=\partial _+W=\partial _+ W'$. Such a splitting is {\it non-trivial} if both $W$ and $W'$ are non-trivial, and {\it doubly-trivial} if both $W$ and $W'$ are trivial. 

\begin{dfn}
A separating surface, $F$, in a 3-manifold, $M$, is {\it strongly irreducible} if every compressing disk for $F$ on one side intersects every compressing disk for $F$ on the other. 
\end{dfn}

\begin{lem}
\label{l:haken}
If $M$ contains a non-trivial or doubly-trivial strongly irreducible Heegaard surface, then $\partial M$ is incompressible in $M$. 
\end{lem}

Lemma \ref{l:haken} follows directly from a Lemma of Haken \cite{haken:68}.

\section{HST splittings}
\label{s:HST}

\begin{dfn} 
A {\it Heegaard-Scharlemann-Thompson (HST) Splitting} of a manifold, $M$, is a sequence of closed, embedded, pairwise disjoint surfaces, $\{F_i\}_{i=0}^n$, such that 
\begin{enumerate}
    \item for each odd $i$ strictly less than $n$, $F_i$ is a non-trivial or doubly-trivial Heegaard splitting of the submanifold of $M$ cobounded by $F_{i-1}$ and $F_{i+1}$. 
    \item if $n$ is odd, then there is a non-trivial compression body, $W$, such that $\partial _+ W=F_n$, and $\partial _- W=F_{n-1}$. 
    \item $\partial M=F_0 \amalg F_n$.
\end{enumerate}
For each odd $i$ the surface $F_i$ is referred to as a {\it thick level}, and for each even $i$, $F_i$ is a {\it thin level}. 
\end{dfn}

\begin{ex}
If $W$ is a compression body, then $\{\partial _-W, \partial _+ W\}$ is an HST splitting of $W$. The only thick level of this HST splitting is $\partial _+W$, and the only thin level is $\partial _- W$. 
\end{ex}

\begin{ex}
If $M$ is a closed 3-manifold, and $F$ is a Heegaard surface for $M$, then an HST splitting of $M$ is $\{\emptyset, F, \emptyset\}$, with a unique thick level, $F$.  
\end{ex}

\begin{ex}
If $M=F \times I$, then $\{F \times \{0\},F \times \{1/2\},F \times \{1\}\}$ is an HST splitting of $M$.
\end{ex}

\begin{dfn}
An HST splitting, $\{F_i\}_{i=0}^n$, is {\it strongly irreducible} if for each odd $i$ strictly less than $n$, the thick level, $F_i$, is strongly irreducible in the submanifold of $M$ cobounded by $F_{i-1}$ and $F_{i+1}$. 
\end{dfn}

\begin{dfn}
\label{d:order}
For any surface, $F$, let $c(F)=\sum \limits _n (2-\chi(F^n))^2$, where $\{F^n\}$ are the components of $F$. If $F_1$ and $F_2$ denote compact, embedded surfaces in a 3-manifold, $M$, then we say $F_1 < F_2$ if $c(F_1) < c(F_2)$.
\end{dfn}

Note that this ordering is defined so that if $F_1$ is obtained from $F_2$ by a compression, then $F_1 < F_2$. 

\begin{dfn}
For any HST splitting, $S=\{F_i\}$, let $c(S)$ denote the ordered set $\{c(F_i)|F_i$ is a thick level$\}$, where repeated integers are included, and the ordering is non-increasing. If $S_1$ and $S_2$ denote HST splittings of some 3-manifold, $M$, then we say $S<S'$ if $c(S)<c(S')$, where the comparison is made lexicographically. 
\end{dfn}

\begin{lem}
Every 3-manifold admits a minimal HST splitting. 
\end{lem}

\begin{proof}
Choose any HST splitting of $M$. Now consider any decreasing sequence of HST splittings of $M$ that begins with the chosen one. This sequence must be finite, since any decreasing sequence of sets of non-negative integers terminates under the lexicographical ordering. 
\end{proof}

The following two Theorems of Scharlemann and Thompson originally appear in \cite{st:94}, with slightly different terminology. One can find proofs of both that are more consistent with the language presented here in \cite{tnbglex}. 

\begin{thm}
\label{t:st}
{\rm [Scharlmann-Thompson \cite{st:94}]} A minimal HST splitting is strongly irreducible.
\end{thm}

\begin{thm}
\label{t:stinceven}
{\rm [Scharlmann-Thompson \cite{st:94}]} Each thin level of a strongly irreducible HST splitting is incompressible. 
\end{thm}

\section{HST splittings of a pair, $(M,K)$}
\label{s:relHST}

To proceed, we must first understand arcs in a compression body. If $W$ is a compression body, recall that $W$ can be built by starting with a product, $F \times I$, and attaching 2- and 3-handles to $F \times \{1\}$. Anything that remains of $F \times \{1\}$ after the attachment becomes part of $\partial _- W$. We say a properly embedded arc, $k$, is {\it straight} in $W$ if $k=\{p\} \times I$, where $p \in F$ is a point such that $\{p\} \times \{1\} \in \partial _-W$. 

We are now ready to generalize the definition of a compression body:

\begin{dfn}
\label{d:Kcompbod}
A {\it $K$-compression body}, $(W;K)$, is 
\begin{enumerate}
    \item A 3-manifold, $W$, which can be obtained by starting with some surface, $F$ (not necessarily connected), forming the product, $F \times I$, attaching some number of 2-handles to $F \times \{1\}$, and capping some (but not necessarily all) remaining 2-sphere boundary components with 3-balls. The boundary component, $F \times \{0\}$, is $\partial _+ W$, and the other boundary component is $\partial _- W$.

    \item  A 1-manifold, $(K, \partial K) \subset (W, \partial W)$, such that
    \begin{enumerate}
        \item $K$ is a disjoint union of embedded arcs
        \item each arc of $K$ has at least one endpoint on $\partial _+W$
        \item if $k$ is an arc of $K$ with $\partial k \subset \partial _+W$, then there is a disk, $D \subset W$, with $\partial D=k \cup \alpha$, $D \cap K=k$, and $D \cap \partial _+ W=\alpha$.  
        \item if $k$ is an arc of $K$ with one endpoint on $\partial _+W$, then $k$ is straight.
        \item For each 2-sphere component, $S$, of $\partial _-W$, $S \cap K \ne \emptyset$.
        \item If $S$ is a 2-sphere component of $\partial _+W$ then either $S$ does not bound a 3-ball component of $W$ or $|K \cap S|\ge 2$.
    \end{enumerate}
\end{enumerate}
A $K$-compression body, $(W;K)$, is {\it non-trivial} if either $W$ is not a product, or at least one arc of $K$ is not straight.
\end{dfn}

\begin{dfn}
If $K$ is a 1-manifold which is properly embedded in a 3-manifold, $M$, then a {\it Heegaard splitting of the pair} $(M,K)$ is an expression of $M$ as a union of $K_i$-compression bodies, $(W_1;K_1)$ and $(W_2;K_2)$, such that $\partial _+W_1=\partial _+ W_2$, and $K=K_1 \cup K_2$. Such a splitting is {\it non-trivial} if both $(W_1;K_1)$ and $(W_2;K_2)$ are non-trivial. 
\end{dfn}

When the context is clear, we will refer to the surface, $\partial _+ W$, as a {\it Heegaard surface} of $(M,K)$.

\begin{dfn} 
If $K$ is a 1-manifold which is properly embedded in a 3-manifold, $M$, then a {\it Heegaard-Scharlemann-Thompson (HST) splitting of the pair} $(M,K)$ is a sequence of closed, embedded, pairwise disjoint surfaces, $\{G_j\}_{j=0}^n$, such that 
\begin{enumerate}
    \item for each odd $j$ between 1 and $n$, there is a non-trivial $K_j$-compression body, $\{(W_j, K_j)\}$, such that $G_j=\partial _+ W_j$, and $G_{j-1}=\partial _- W_j$.
    \item for each even $j$ between 2 and $n$, there is a non-trivial $L_j$-compression body, $\{V_j, L_j\}$, such that $G_{j-1}=\partial _+ V_j$, and $G_j=\partial _- V_j$.
    \item $\partial M=G_0 \amalg G_n$
    \item $K=\bigcup \limits _j K_j \cup L_j$.
\end{enumerate}
For each odd $j$ the surface $G_j$ is referred to as a {\it thick level}, and for each even $j$, $G_j$ is a {\it thin level}. 
\end{dfn}

\begin{ex}
\label{e:thin}
Suppose that $K \subset S^3$ is an arbitrary knot or link with no trivial components, and $h$ is some standard height function on $S^3$ (so that for each $p \in (0,1)$, $h^{-1}(p)$ is a 2-sphere), which is a Morse function when restricted to $K$. Let $\{ q'_j \}$ denote the critical values of $h$ restricted to $K$, and let $q_j$ be some point in the interval $(q'_j, q'_{j+1})$. The following terminology is standard in {\it thin position} arguments (see \cite{gabai:87}).

If $j$ is such that $|K \cap h^{-1}(q_j)|>|K \cap h^{-1}(q_{j-1})|$ and $|K \cap h^{-1}(q_j)|>|K \cap h^{-1}(q_{j+1})|$, then we say the surface $h^{-1}(q_j)$ is a {\it thick level} of $K$. Similarly, if $|K \cap h^{-1}(q_j)|<|K \cap h^{-1}(q_{j-1})|$ and $|K \cap h^{-1}(q_j)|<|K \cap h^{-1}(q_{j+1})|$, then we say the surface $h^{-1}(q_j)$ is a {\it thin level} of $K$. 

Suppose there are $n$ thick levels for $K$. Let $G_0=G_{2n}=\emptyset$, $\{G_{2j-1}\}_{j=1}^n$ denote the set of thick levels of $K$, and $\{G_{2j}\}_{j=0}^{n-1}$ denote the set of thin levels. Then $\{G_j\}_{j=0}^{2n}$ is an HST splitting of $(S^3,K)$, and the thick and thin levels of this HST splitting are precisely the thick and thin levels of $K$.  
\end{ex}

\begin{dfn}
If $(K,\partial K) \subset (M,\partial M)$ is some properly embedded 1-manifold then let $M_K$ denote $M$ with a regular neighborhood of $K$ removed. If $X$ is any subset of $M$ then let $X_K=M_K \cap X$.
\end{dfn}

\begin{dfn}
Suppose $F$ is a surface in a 3-manifold $M$, and $K$ is some properly embedded 1-manifold in $M$, transverse to $F$. A {\it relative compressing disk} for $F$ is a disk, $D$, such that $\partial D=\alpha \cup \beta$, where $D \cap K=\alpha$, and $D \cap F =\beta$. If there are no relative compressions for some surface, then we say that it is {\it relatively incompressible}. To {\it perform a relative compression} is to use a relative compressing disk to guide an isotopy, thereby reducing the number of intersections with $K$ by exactly two. 
\end{dfn}

\begin{dfn}
\label{d:weakreduce}
Suppose $F$ is a separating surface in a 3-manifold, $M$, and $K$ is some properly embedded 1-manifold in $M$, transverse to $F$. Then $F$ is {\it weakly reducible with respect to $K$} if there exist disks, $D$ and $E$, on opposite sides of $F$, such that either
    \begin{enumerate}
        \item $D$ and $E$ are disjoint compressing disks for $F_K$ in $M_K$
        \item $D$ and $E$ are relative compressing disks for $F$ that are either disjoint, or meet in a point of $K$
        \item $D$ is a relative compressing disk for $F$, and $E$ is a disjoint compressing disk for $F_K$ in $M_K$
        \item $D$ is a compressing disk for $F_K$ in $M_K$, and $E$ is a disjoint relative compressing disk for $F$. 
    \end{enumerate}
We say $F$ is {\it strongly irreducible with respect to $K$} if it is not weakly reducible with respect to $K$.
\end{dfn}

\begin{ex}
\label{e:thin3}
Let $K$, $h$, and $\{q_j\}$ be as in Example \ref{e:thin}. The {\it width} of $K$ is defined to be the quantity $\sum \limits _j |K \cap h^{-1}(q_j)|$. A knot is said to be in {\it thin position} if $h$ is chosen so that the width of $K$ is minimal (see \cite{gabai:87}).

Recall from Example \ref{e:thin} that any choice of $h$ induces an HST splitting, $\{G_j\}$, of $(S^3,K)$. If, when $h$ is chosen so as to minimize the width of $K$, it turns out that $\{G_j\}=\{\emptyset, G, \emptyset\}$ ({\it i.e.} $K$ is also in {\it bridge position} with respect to $h$), then we claim that $G$ is a strongly irreducible Heegaard splitting of $(S^3,K)$. 

Suppose this were not the case. Then $G$ is weakly reducible with respect to $K$. Let $D$ and $E$ be a pair of disks, as in the above definition of weak reducibility. We now examine each of the possible cases in that definition:

    \begin{enumerate}
        \item {\it $D$ and $E$ are relative compressing disks for $G$ that are either disjoint, or meet in a point of $K$.} Then there is an isotopy of $K$ which reduces its width, as depicted in Figure \ref{f:thinmove}.

        \begin{figure}[htbp]
        \psfrag{h}{$D$}
        \psfrag{l}{$E$}
        \vspace{0 in}
        \begin{center}
        \epsfxsize=2 in
        \epsfbox{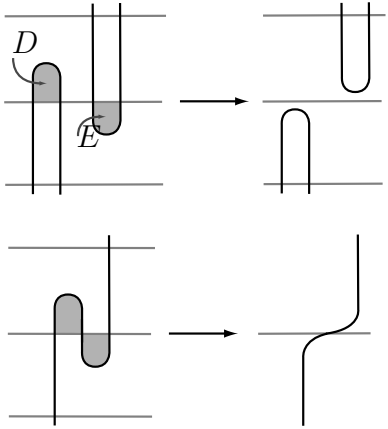}
        \caption{Moves which produce lower width.}
        \label{f:thinmove}
        \end{center}
        \end{figure}

        \item {\it $D$ is a relative compressing disk for $G$, and $E$ is a disjoint compressing disk for $G_K$ in $M_K$.} Then $\partial E$ bounds disks, $E'$ and $E''$, on $G$ (since $G$ is a sphere). One of these disks ($E'$, say) is disjoint from $D$. As $E$ is a compressing disk for $G_K$, it must be that there is a relative compressing disk inside the sphere bounded by $E \cup E'$. However, this relative compressing disk is disjoint from $D$, so we are reduced to Case 1.
        \item {\it $D$ is a compressing disk for $G_K$ in $M_K$, and $E$ is a disjoint relative compressing disk for $F$.} This case and the previous one are completely symmetric. 
        \item {\it $D$ and $E$ are disjoint compressing disks for $G_K$ in $M_K$.} As in Case 2, $\partial E$ bounds disks, $E'$ and $E''$ on $G$. Also, $\partial D$ bounds disks, $D'$ and $D''$ on $G$. Since $\partial D \cap \partial E=\emptyset$, one of the disks, $E'$ or $E''$ is disjoint from one of the disks, $D'$ or $D''$. Say $E' \cap D'=\emptyset$. As in Case 2, there is a relative compression, $L$, inside the ball bounded by $E \cup E'$, and a relative compression, $H$, inside the ball bounded by $D \cup D'$. But $H \cap L=\emptyset$, so we are again reduced to Case 1 above.
    \end{enumerate}
\end{ex}

In Lemma \ref{l:untang} we will present a generalization of this example for graphs in 3-manifolds, as opposed to 1-manifolds in $S^3$.

Proofs of the next two Lemmas originally appeared in \cite{bachman:98} with different terminology. One can find proofs of both that are more consistent with the language presented here in \cite{tnbglex}. 

\begin{lem}
\label{l:essential}
If $F$ is a strongly irreducible Heegaard surface of $(M,K)$ and $(S,\partial S) \subset (M, \partial M \cup K)$ is an essential surface then $F$ can be isotoped rel $K$ so that no arc of $S \cap F$ cobounds a relative compressing disk for $F$, and so that no loop of $S \cap F$ bounds a compressing disk for $F_K$ in $M_K$.
\end{lem}

\begin{lem}
\label{l:relhaken}
If $F$ is a strongly irreducible Heegaard surface of $(M,K)$ then $(\partial M)_K$ is incompressible in $M_K$, and $\partial M$ is relatively incompressible in $M$.
\end{lem} 

\begin{dfn}
An HST splitting, $\{G_j\}_{j=0}^n$ of a pair, $(M,K)$, is {\it strongly irreducible} if for each odd $j$ strictly less than $n$, the thick level $G_j$ is strongly irreducible with respect to $K$ in the submanifold of $M$ cobounded by $G_{j-1}$ and $G_{j+1}$. 
\end{dfn}

We now present one further generalization of the concept of a HST splitting. If $\Gamma$ is any graph then let $\Gamma ^0$ and $\Gamma ^1$ denote the sets of vertices and edges of $\Gamma$. If we are discussing a particular graph, $\Gamma$, embedded in a 3-manifold, $M$, then $M^*$ will denote $M$ with a regular neighborhood of $\Gamma ^0$ removed. 

\begin{dfn}
Let $\Gamma$ be a properly embedded graph in a 3-manifold, $M$. We define an {\it HST splitting of $(M,\Gamma)$} to be an HST splitting of $(M^*, \Gamma ^1)$.
\end{dfn}

\section{Underlying HST Splittings}

If $\Gamma$ is a graph properly embedded in an irreducible 3-manifold, $M$, and we are given an HST splitting of $(M,\Gamma)$, then we can define an HST splitting of $M$ by ``forgetting" $\Gamma$. 

\begin{dfn}
\label{d:underlyingHST}
If $\{G_j\}_{j=0}^n$ is an HST splitting of $(M,\Gamma)$, where $M$ is an irreducible 3-manifold other than $B^3$ or $S^3$, then we define its {\it underlying HST splitting}, $[\{G_j\}]$, as follows: 
\begin{enumerate}
    \item Let $s:\{0,...,n\} \rightarrow \{0,...,m\}$ be the onto, monotone function, such that $s(i)=s(j)$ iff the submanifold of $M$ co-bounded by the non-$S^2$ components of $G_i$ and $G_j$ is a product. 
    \item For each $i$ between 0 and $m$, choose some $j \in s^{-1}(i)$, and let $F_i$ denote the non-$S^2$ components of $G_j$.  
    \item Let $\sigma$ be the maximal subset of $\{0,...,m\}$ such that $\{F_i\}_{i \in \sigma}$ is an HST splitting of $M$. 
    \item Define $[\{G_j\}]=\{F_i\}_{i \in \sigma}$.
\end{enumerate}
We leave it to the reader to check that $[\{G_j\}]$ is well defined up to isotopy.
\end{dfn}

\begin{dfn}
Let $\Gamma$ be a properly embedded graph in a 3-manifold, $M$. For any surface, $F \subset M$, let $c(F;\Gamma)=\sum \limits _n \left(2-\chi(F^n_{\Gamma})\right)^2$, where $\{F^n\}$ are the components of $F$, and $F^n_{\Gamma}$ denotes $F^n$ with a neighborhood of $\Gamma$ removed. If $F_1$ and $F_2$ denote compact, embedded surfaces in $M$, 
then we say $F_1 <_{\Gamma} F_2$ if $c(F_1;\Gamma) < c(F_2;\Gamma)$.
\end{dfn}

In the next Lemma we relate the thick and thin levels of an HST splitting of a pair to those of its underlying HST splitting. This will be an important step in the eventual production of almost normal surfaces.

\begin{lem}
\label{l:underlyingodd}
Let $\Gamma$ be a properly embedded graph in an irreducible 3-manifold $M$ other than $S^3$ or $B^3$. Let $\{G_j\}$ be an HST splitting of $(M,\Gamma)$, and suppose $\{F_i\}=[\{G_j\}]$. Then for each thick (thin) level, $F_p$, of $\{F_i\}$ there is a thick (thin) level of $\{G_j\}$ which becomes parallel to $F_p$ when all $S^2$ components are removed. 
\end{lem}

\begin{proof}
Let $s$ and $\sigma$ be as given in Definition \ref{d:underlyingHST}. Let $x$ be the $k^{\rm th}$ element of $\sigma$, for some odd number, $k$. $F_x$ is parallel to the non-$S^2$ components of $G_j$, for all $j \in s^{-1}(x)$. Since $s$ is monotone, $s^{-1}(x)$ either contains an odd number, or contains only one element. Assume the latter is true, and let $y=s^{-1}(x)$. The assumption that $s^{-1}(x)$ has only one element implies the non-$S^2$ components of $G_{y-1}$, $G_y$, and $G_{y+1}$ (which we denote $G_{y-1}'$, $G_y'$, and $G_{y+1}'$) are non-parallel. 

Now, let $x^-$ denote the predecessor of $x$ in $\sigma$. Since $k$ is odd, $F_x>F_{x^-}$. Since $y$ is even, $G_y <_{\Gamma} G_{y-1}$. As $y-1 \notin s^{-1}(x)$, we know that $G_y'$ and $G_{y-1}'$ are not parallel, so it must be that $G_y' < G_{y-1}'$. But $y-1 \notin s^{-1}(x)$ also implies that $F_{x-1}$ is parallel to $G_{y-1}'$. Since $x^-$ is the predecessor of $x$ in $\sigma$, $x^- \le x-1 \le x$. Since $\sigma$ is maximal, $F_x>F_{x^-}$ implies $F_x>F_{x-1} \ge F_{x^-}$. This is a contradiction, as $F_x=G_y'$, $F_{x-1}$ is parallel to $G_{y-1}'$, and $G_y' < G_{y-1}'$.

The proof in the case that $k$ is even is completely analogous. 
\end{proof}

\section{Strongly irreducible HST Splittings of pairs}

In this section we examine the properties of strongly irreducible HST splittings of pairs and prove an important Theorem about their existence. We begin with a result which is completely analogous to Theorem \ref{t:stinceven}. Its proof is similar to that of Theorem 7.4. of \cite{tnbglex}.

\begin{thm}
\label{t:minwrtk}
Let $K$ be a properly embedded 1-manifold in a 3-manifold, $M$. Each thin level of a strongly irreducible HST splitting of $(M,K)$ is incompressible in $M_K$, and relatively incompressible in $M$. 
\end{thm}

\begin{proof}
Let $\{G_j\}_{j=0}^m$ be a strongly irreducible HST splitting of $(M,K)$. For each odd $j$, let $M_j$ denote the submanifold of $M$ cobounded by $G_{j-1}$ and $G_{j+1}$. Suppose $c$ is a loop that bounds a compressing disk, or an arc which cobounds a relative compressing disk, for some surface, $G_q$, where $q$ is even. If $c$ is a loop, then let $C$ denote a disk in $M$, such that $\partial C=c$. If $c$ is an arc, then let $C$ be a disk in $M$ such that $\partial C=\gamma \cup c$, where $K \cap C=\gamma$. In either case, choose $C$ so that $|C \cap (\bigcup \limits _{{\rm even} \ j} G_j)|$ is minimal. 

It follows from Lemma \ref{l:relhaken} that for each odd $j$, $\partial M_j$ is incompressible and relatively incompressible in $M_j$. Hence, $C$ cannot lie entirely in $M_{q-1}$ or $M_{q+1}$. We conclude then that there is some loop or arc of intersection of the interior of $C$ with $\bigcup \limits _{{\rm even} \ j} G_j$. Let $\alpha$ denote an innermost such loop. Let $C'$ denote the subdisk of $C$ bounded by $\alpha$. $C'$ lies in $M_p$, for some odd number, $p$. As $\partial M_p$ is incompressible in $M_p$, $\alpha$ must bound a disk, $A$, on $\partial M_p$. 

Now, let $\beta$ be an innermost loop of $C \cap A$, and let $A'$ be the subdisk of $A$ bounded by $\beta$. Then we can use $A'$ to surger $C$, and thereby obtain a new disk, with the same boundary as $C$, contradicting our minimality assumption.

We conclude then that $c$ cannot be a loop, and if $c$ is an arc, then $C$ contains no loops of intersection with $\bigcup \limits _{{\rm even} \ j} G_j$. Let $\delta$ then denote an arc of intersection which is outermost on $C$. $\delta$ and a subarc of $\alpha$ cobound a subdisk, $C''$, of $C$. $C''$ is then a relative compressing disk for $\partial M_j$, for some $j$, contradicting Lemma \ref{l:relhaken}.
\end{proof}

\begin{cor}
If $\Gamma$ is a properly embedded graph in an irreducible 3-manifold, $M$, then each thin level of a strongly irreducible HST splitting of $(M,\Gamma)$ is incompressible in $M_{\Gamma}=M-Nbhd(\Gamma)$, and relatively incompressible in $M^*$. 
\end{cor}

\begin{dfn}
For any HST splitting, $S=\{F_i\}$, of $(M,\Gamma)$, let $c(S;\Gamma)$ denote the ordered set $\{c(F_i;\Gamma)|F_i$ is a thick level of $S\}$, where repeated integers are included, and the ordering is non-increasing. If $S_1$ and $S_2$ denote HST splittings of $(M,\Gamma)$, then we say $S_1 <_{\Gamma} S_2$ if $c(S_1;\Gamma)<c(S_2;\Gamma)$, where the comparison is made lexicographically. 
\end{dfn}

\begin{lem}
\label{l:untang}
Let $\Gamma$ be a properly embedded graph in an irreducible 3-manifold, $M$. If $\{F_i\}$ is a strongly irreducible HST splitting of $M$ then there exists a strongly irreducible HST splitting, $\{G_j\}$, of $(M,\Gamma)$ such that $[\{G_j\}]=\{F_i\}$.
\end{lem}

\begin{proof}
Among all HST splittings of $(M,\Gamma)$ whose underlying HST splitting is $\{F_i\}$, choose one, $\{G_j\}_{j=0}^m$, that is smallest. We would like to show that $\{G_j\}$ is strongly irreducible. For each odd $j$ strictly less than $m$, let $M_j$ denote the submanifold of $M$ cobounded by $G_{j-1}$ and $G_{j+1}$. Let $G_p$ be a thick level, so that $p$ is odd. We need to show that $G_p$ is strongly irreducible with respect to $\Gamma \cap M_p$ in $M_p$. If not, then there are disks, $D$ and $E$, in $M_p$, as in Definition \ref{d:weakreduce}. 

Note that we may assume that $\partial D$ does not bound a disk, $D'$, on $G_p$ such that $|D \cup D'|$ bounds a ball containing a single unknotted arc of $\Gamma ^1$. If this were the case then there is a relative compressing disk, $D^*$, inside the ball bounded by $D \cap D'$. As $D^* \cap E=\emptyset$, we may use $D^*$ instead of $D$. We make a similar choice for $E$. Our reason for such choices for $D$ and $E$ will be made clear shortly. 

Let $G_D$, $G_E$, and $G_{DE}$ be the surfaces obtained from $G_p$ by (relative) compression along $D$, (relative) compression along $E$, and (relative) compression along both $D$ and $E$. By our choice of $D$ from the previous paragraph we know that if some component $S$ of $G_D$ bounds a ball in $M_p$ then $|K \cap \Gamma ^1|\ge 2$. This a necessary condition for $G_D$ and $G_{p-1}$ to cobound a $\Gamma ^1$-compression body. (See condition 2(f) of Definition \ref{d:Kcompbod}.)

We now produce a new HST splitting from $\{G_j\}$. There are four cases:

\begin{enumerate}
    \item $G_D \ne G_{p-1}$, $G_E \ne G_{p+1}$. 
    \\ Remove $G_p$ from $\{G_j\}$. In its place, insert $\{G_D, G_{DE}, G_E\}$ and reindex.  
    \item $G_D = G_{p-1}$, $G_E \ne G_{p+1}$. 
    \\ Replace $\{G_{p-1}, G_p\}$ with $\{G_{DE},G_E\}$.
    \item $G_D \ne G_{p-1}$, $G_E = G_{p+1}$. 
    \\ Replace $\{G_p, G_{p+1}\}$ with $\{G_D,G_{DE}\}$.
    \item $G_D = G_{p-1}$, $G_E = G_{p+1}$. 
    \\ Replace $\{G_{p-1}, G_p, G_{p+1}\}$ with $G_{DE}$ and reindex. 
\end{enumerate}

In all cases it is a routine matter to check that we have defined a new HST splitting of $(M,\Gamma)$ which is smaller than $\{G_j\}$. What remains to be checked is that the underlying HST splitting is still $\{F_i\}$. We claim that the only way in which this can fail is if the non-$S^2$ components of both $G_D$ and $G_E$ are not parallel to the non-$S_2$ components of $G_p$. The proof of this assertion is similar in all four cases. We present a proof for the first (and most difficult) case and leave the remaining cases as exercises. 

Let $\sigma$ be as in Definition \ref{d:underlyingHST}. If the non-$S^2$ components of two surfaces, $A$ and $B$, are parallel then we write $A \sim B$. Suppose that $G_p \sim G_D$. Then $D$ is either a relative compressing disk, or $\partial D$ bounds a disk on $G_p$ that meets $\Gamma$ (and hence $\partial D$ is essential on $(G_p)_{\Gamma}$, but inessential on $G_p$). Since $G_{DE}$ is obtained from $G_E$ by (relative) compression along $D$, we conclude that $G_E \sim G_{DE}$. 

In Case 1 above, the new HST splitting that we construct looks like
\[\{...G_{p-1}, G_D, G_{DE}, G_E, G_{p+1},...\}\]
As $G_D \sim G_p$, and $G_{DE} \sim G_E$, the underlying HST splitting of this is the same as the underlying HST splitting of the following sequence:
\[\{...G_{p-1}, G_p, G_E, G_{p+1},...\}\]
Since $E$ is a disk in the submanifold of $M$ cobounded by $G_p$ and $G_{p+1}$ we will not include the index of $G_E$ in $\sigma$ when forming the underlying HST splitting of this sequence. Hence, the underlying HST splitting of this last sequence is the same as that of $\{G_j\}$. 

What we have established is that if the underlying HST splitting of the new HST splitting of $(M, \Gamma)$ constructed above is not equal to $\{F_i\}$ then both $D$ and $E$ were compressing disks for $G_p$. As $G_p$ is parallel to some thick level of $\{F_i\}$ (by Lemma \ref{l:underlyingodd}), and $D$ and $E$ are disjoint, this contradicts our assumption that $\{F_i\}$ was strongly irreducible.
\end{proof}

\section{Normal Surfaces: Definitions}
\label{s:normal}

\begin{figure}[htbp]
        \vspace{0 in}
        \begin{center}
        \epsfxsize=3.25 in
        \epsfbox{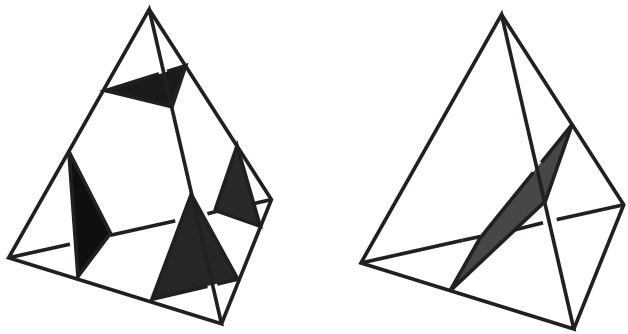}
        \caption{Normal Disks.}
        \label{f:Normal}
        \end{center}
        \end{figure}

In this section, we discuss the necessary background material on normal surfaces. A {\it normal curve} on the boundary of a tetrahedron is a simple loop which is transverse to the 1-skeleton, made up of arcs which connect distinct edges of the 1-skeleton. The {\it length} of such a curve is simply the number of times it crosses the 1-skeleton. A {\it normal disk} in a tetrahedron is any embedded disk, whose boundary is a normal curve of length three or four, and whose interior is contained in the interior of the tetrahedron, as in Figure \ref{f:Normal}.

A {\it normal surface} in $M$ is the image of an embedding, $p$, of some closed surface, $F$, into $M$, such that $p(F)$ is a union of normal disks. In addition, we say $p(F)$ is an {\it almost normal surface} if it consists of all normal disks, plus one additional piece in one tetrahedron. This piece can be either a disk with normal boundary of length 8 (depicted in Figure \ref{f:Octagon}), or two normal disks connected by a single unknotted tube (as in Figure \ref{f:1Tube}). Almost normal surfaces were first explored by Rubinstein in \cite{rubinstein:93}, and later used by Thompson \cite{thompson:94} and Stocking \cite{stocking:96}. In \cite{bachman:98}, we produce new applications of almost normal surfaces by generalizing them to include surfaces with non-empty boundary.

\section{Almost Normal Surfaces and HST splittings of $(M,T^1)$}
\label{s:nlan}

One application of the results we have discussed thus far comes about when we let $\Gamma$ be the 1-skeleton of a triangulation of a 3-manifold, $M$. If $T$ is such a triangulation then we denote the $n$-skeleton of $T$ by $T^n$. We now focus on strongly irreducible HST splittings of $(M,T^1)$.

\begin{dfn}
\label{bubble}
Suppose $F$ is an embedded surface in $M$, equipped with a triangulation, $T$. A {\it bubble} for $F$ is a ball, $B$, such that $\partial B= D_1 \cup D_2$, where $D_1$ and $D_2$ are disks, $D_1$ is contained in a single tetrahedron, $F \cap B=D_2$, $D_2 \cap T^2 \ne \emptyset$, and $D_2 \cap T^1 = \emptyset$.
\end{dfn}

\begin{lem}
\label{l:nobubbles}
Suppose $M$ is a submanifold of a 3-manifold, $M'$, with triangulation, $T$, which is bounded by a normal or almost normal surface. Let $\{G_j\}$ be a strongly irreducible HST splitting of $(M,T^1)$. Then we may isotope $\{G_j\}$ so that for all $j$, $G_j$ does not contain any bubbles. 
\end{lem}

\begin{proof}
Suppose $B$ is a bubble for $G_j$, for some $j$, where $\partial B=D_1 \cup D_2$, as in Definition \ref{bubble}. We can use $B$ to guide an isotopy from $D_2$ to $D_1$. This may push other surfaces of $\{G_j\}$ which had non-empty intersection with $int(B)$, but it can only destroy bubbles for those surfaces, too. The isotopy is supported on a neighborhood of $B$, which is disjoint from $T^1$. Hence, if $\{G_j\}$ was strongly irreducible before the isotopy then it will remain so. Since there are a finite number of surfaces in $\{G_j\}$, and a finite number of bubbles for each, we eventually arrive at an HST splitting with the desired properties.
\end{proof}

\begin{thm}
\label{t:normal}
Suppose $M$ is a submanifold of a 3-manifold, $M'$, with triangulation, $T$, which is bounded by a normal or almost normal surface. Let $\{G_j\}$ be a strongly irreducible HST splitting of $(M,T^1)$. Then we may isotope $\{G_j\}$ so that each thin level is a normal surface.
\end{thm}

\begin{proof}
We begin by using Lemma \ref{l:nobubbles} to isotope $\{G_j\}$, removing all bubbles.

Now, let $j$ be some even number, let $\tau$ be some tetrahedron in $T$, and let $\Delta$ be a face of $\tau$. First, we examine the possibilities for $G_j \cap \Delta$. Let $\gamma$ be an innermost simple closed curve, bounding a disk, $D_1$ in $\Delta$. By Theorem \ref{t:minwrtk}, $\gamma$ must bound a disk, $D_2$, in $G_j$. $M_{T^1}$ is irreducible (it's a handlebody), so $D_1 \cup D_2$ bounds a bubble for $G_j$. This is a contradiction, so we see no simple closed curves in any face.

If there are any curves which run from one edge of $\Delta$ to itself, then there is an outermost such one. Let $D$ denote the sub-disk it cuts off in $\Delta$. Then $D$ is a relative compressing disk for $G_j$, also contradicting Theorem \ref{t:minwrtk}. We conclude that $G_j \cap \Delta$ is a collection of normal arcs.

We now consider the possibilities for $G_j \cap \partial \tau$. It is easy to show that the only possibilities for normal loops are curves of length 3, or 4n (see, for example, \cite{thompson:94}). If there are any curves of length greater than 4, then there must be a disk, $D$, such that $\partial D=\alpha \cup \beta$, where $D \cap T^1 =\alpha$, and $D \cap G_j =\beta$ (see \cite{thompson:94}). This is a relative compressing disk for $G_j$, which is again a contradiction. We conclude that $G_j \cap \partial \tau$ consists of normal loops of length 3 and 4.

Finally, it follows from Theorem \ref{t:minwrtk} that every loop of $G_j \cap \partial \tau$ bounds a disk on $G_j$. Since we have already ruled out simple closed curves in faces of $\tau$, such disks must lie entirely inside $\tau$. We conclude $G_j$ is a normal surface.
\end{proof}

Our goal now is to show that once bubbles are removed from the thick levels of $\{G_j\}$, they become almost normal in $M$. First, we shall need a few lemmas.

\begin{lem}
\label{l:nlarcs}
Suppose $M$ is a submanifold of a 3-manifold, $M'$, with triangulation, $T$, which is bounded by a normal or almost normal surface. Let $\{G_j\}_{j=0}^m$ be a strongly irreducible HST splitting of $(M,T^1)$, where $\partial M$ normal implies $m$ is even, and $\partial M$ almost normal implies $m$ is odd. Then we may isotope $\{G_j\}$, so that each thin level is a normal surface, and so that for each thick level, $G_j$, $G_j \cap \Delta$ is a collection of normal arcs, for every 2-simplex, $\Delta$, of $T$. 
\end{lem}

\begin{proof}
We begin by appealing to Lemma \ref{l:nobubbles} and Theorem \ref{t:normal} to isotope $\{G_j\}$ so that there are no bubbles for any of its surfaces, and so that each thin level is normal. We now produce a further isotopy of $\{G_j\}$, which fixes all thin levels, and which results in an HST splitting with all of the desired properties. The isotopy will be defined individually for each odd value of $j$ ({\it i.e.}, for each thick level). 

Let $j$ be some odd number, and let $M_j$ denote the submanifold of $M$ cobounded by $G_{j-1}$ and $G_{j+1}$. The strong irreducibility of $\{G_j\}$ implies that $G_j$ is a strongly irreducible Heegaard surface for $(M_j, T^1)$. Lemma \ref{l:essential} now implies that a standard innermost disk/outermost arc argument can be used to isotope $G_j$ so that no arc of $T^2 \cap G_j$ cobounds a relative compressing disk for $G_j$ in $M_j$, and so that no loop of $T^2 \cap G_j$ bounds a compressing disk for $(G_j)_{T^1}$ in $(M_j)_{T^1}$. 

Now, suppose there are non-normal components of $G_j \cap \Delta$, for some 2-simplex, $\Delta$. If there are loops of $G_j \cap \Delta$, then let $\gamma$ be a loop which is innermost on $\Delta$. Then $\gamma$ bounds a subdisk, $D$, of $\Delta$. But, as no loop of $T^2 \cap G_j$ bounds a compressing disk for $(G_j)_{T^1}$ in $(M_j)_{T^1}$, it must be that $\gamma$ also bounds a disk, $D'$, on $(G_j)_{T^1}$. But then $D \cup D'$ bounds a bubble for $G_j$. As all such bubbles have been removed, this is a contradiction. 

We conclude then that if there are non-normal components of $G_j \cap \Delta$, then they must be arcs which run from some edge of $T^1$ back to itself. But then if $\alpha$ is such an arc which is outermost on $\Delta$, then $\alpha$ cobounds a subdisk of $\Delta$ which is a relative compressing disk for $G_j$, a contradiction. 
\end{proof}

\begin{lem}
\label{348}
Let $M,T$ and $\{G_j\}$ be as given by the conclusion of Lemma \ref{l:nlarcs}. Then for each odd $j$, $G_j$ meets the boundary of every tetrahedron in normal curves of length 3, 4, and at most one curve on at most one tetrahedron of length 8.
\end{lem}

This lemma is taken straight from \cite{thompson:94}. We refer the reader to this paper for its proof. The necessary assumptions are that $G_j$ meets every tetrahedron in normal arcs, and that there is no pair of disks which form a weak reduction for $G_j$.

A connected surface in a ball, $B$, is {\it unknotted} if it is isotopic to a neighborhood of the graph obtained by coning $n$ points on $\partial B$ to a point in the interior of $B$. For the proof of Theorem \ref{t:an} we will need the following result:

\begin{thm}
\label{t:sch}
{\rm [Scharlemann \cite{scharlemann:97}]} Suppose $F$ is a strongly irreducible Heegaard surface for a 3-manifold, $M$, and $B$ is a ball in $M$. If $\partial B \backslash F$ is incompressible in $M \backslash F$, then $F \cap B$ is connected and unknotted. 
\end{thm}

To prove Theorem \ref{t:an}, we will need to restate the previous result for Heegaard surfaces of $(M,K)$:

\begin{thm}
\label{t:relsch}
Suppose $F$ is a strongly irreducible Heegaard surface of $(M;K)$ and $B$ is a ball in $M_K$. If $\partial B \backslash F$ is incompressible in $M_K \backslash F$, then $F \cap B$ is connected and unknotted. 
\end{thm}

The proof of this Theorem is exactly the same as Scharlemann's proof of Theorem \ref{t:sch}. The point is that there is very little difference between a Heegaard surface for $M$ and a Heegaard surface for $(M,K)$, as viewed from inside a ball that lies entirely in the complement of $K$.

\begin{thm}
\label{t:an}
Suppose $M$ is a submanifold of a 3-manifold, $M'$, with triangulation, $T$, which is bounded by a normal or almost normal surface. Let $\{G_j\}_{j=0}^m$ be a strongly irreducible HST splitting of $(M,T^1)$, where $\partial M$ normal implies $m$ is even, and $\partial M$ almost normal implies $m$ is odd. Then we may isotope $\{G_j\}$, so that each thin level is a normal surface and each thick level is almost normal.
\end{thm}

\begin{proof}
Isotope $\{G_j\}$ to satisfy the conclusion of Lemma \ref{l:nlarcs}. Let $j$ be some odd number, so that $G_j$ is a thick level. Let $\tau$ be some tetrahedron in $T$. Let $S$ be a copy of $\partial \tau$, pushed slightly into $\tau$. We first claim that $S \backslash G_j$ is incompressible in $M_{T^1}$, to the outside of $S$ (the side which is not contained in $\tau$). Suppose not, and let $D$ denote a compressing disk which meets $\partial \tau$ minimally. As $D$ lies in $M_{T^1}$, it must be that $D \cap \Delta$ is a (possibly empty) collection of simple closed curves, for each 2-simplex, $\Delta$, of $\partial \tau$. Let $\gamma$ denote an innermost such simple closed curve. Then $\gamma$ bounds a subdisk of $\Delta$ which is disjoint from $G_j$ (since $G_j$ meets $\Delta$ only in normal arcs). Hence, we can use this subdisk to surger $D$, thereby obtaining a compressing disk for $S \backslash G_j$ which meets $\partial \tau$ in one fewer curve. Our conclusion, then, is that $D \cap \partial \tau=\emptyset$. But then $D$ lies between $S$ and $\partial \tau$. As these two surfaces are parallel, this is impossible. 

Now, choose a complete collection of compressing disks for $S \backslash G_j$ inside $\tau \backslash G_j$, and surger $S$ along this collection. We obtain in this way a collection of spheres, $\{ S_1, ... ,S_n \}$. $S_i$ bounds a ball, $B_i$, in $\tau$, and by definition, $\partial B_i \backslash G_j$ is incompressible in the complement of $(G_j) _ {T^1}$ in $M _ {T^1}$. These are the conditions necessary to apply Theorem \ref{t:relsch}. The conclusion is that inside each $B_i$, $G_j$ is a connected surface, which looks like the neighborhood of a graph which is the cone on some collection of points in $\partial B_i$. So, in particular, if $G_j \cap \partial B_i$ is a single curve, then it bounds a disk in $B_i$, and hence so does the corresponding curve in $\partial \tau$.

        \begin{figure}[htbp]
        \vspace{0 in}
        \begin{center}
        \epsfxsize=5 in
        \epsfbox{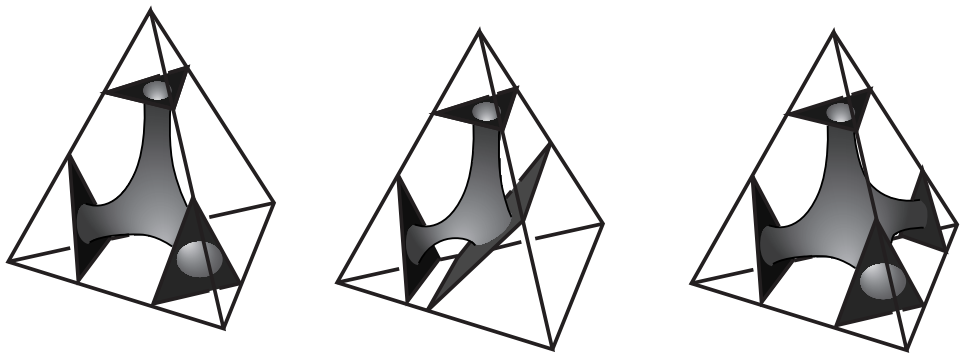}
        \caption{Possibilities when $G_j \cap \partial B_i$ consists of 3 or more curves.}
        \label{f:2Tubes}
        \end{center}
        \end{figure}

Suppose there is some $i$ such that $G_j \cap \partial B_i$ consists of three or more curves, of length 3 or 4. The only ways this can happen are shown in figure \ref{f:2Tubes}. In all cases we see a compressing disk on one side of $G_j$ which is disjoint from a relative compressing disk on the other side (see figure \ref{f:HLComp}). 

        \begin{figure}[htbp]
        \vspace{0 in}
        \begin{center}
        \epsfxsize=3.25 in
        \epsfbox{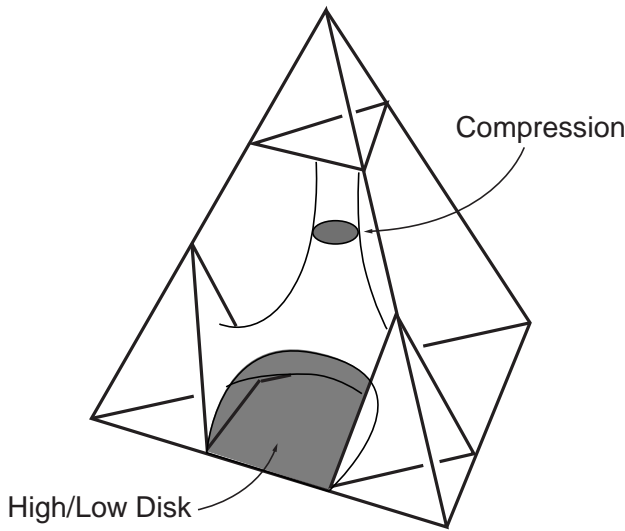}
        \caption{A disjoint compression and relative compression.}
        \label{f:HLComp}
        \end{center}
        \end{figure}

Now suppose that for some $i$, $G_j \cap \partial B_i$ consists of two normal curves, of length 3 or 4. Theorem \ref{t:relsch} tells us that the picture must be two normal disks, tubed together by a single unknotted tube, as in figure \ref{f:1Tube}. Note that in this situation, we see a relative compressing disk on one side, and a compressing disk on the other. Hence, there cannot be more than one place where we see this picture. Otherwise, we'd see either two disjoint compressing disks on opposite sides, or a compressing disk on one side disjoint from a relative compressing disk on the other. Neither of these situations can happen for a surface which is strongly irreducible with respect to $T^1$.

        \begin{figure}[htbp]
        \vspace{0 in}
        \begin{center}
        \epsfxsize=3.75 in
        \epsfbox{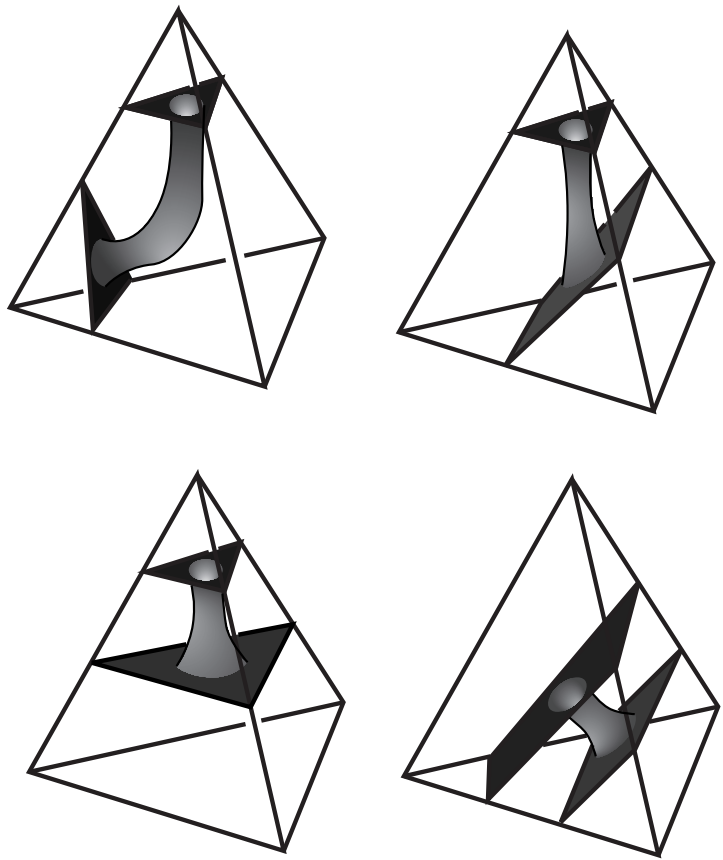}
        \caption{Possibilities when $G_j \cap \partial B_i$ consists of 2 curves.}
        \label{f:1Tube}
        \end{center}
        \end{figure}

Furthermore, suppose $G_j \cap \partial \tau$ contains a curve of length 8. Then we see relative compressing disks on both sides as in figure \ref{f:Octagon}, and hence, there cannot be a tube anywhere else (including attached to this disk!).

        \begin{figure}[htbp]
        \vspace{0 in}
        \begin{center}
        \epsfxsize=3 in
        \epsfbox{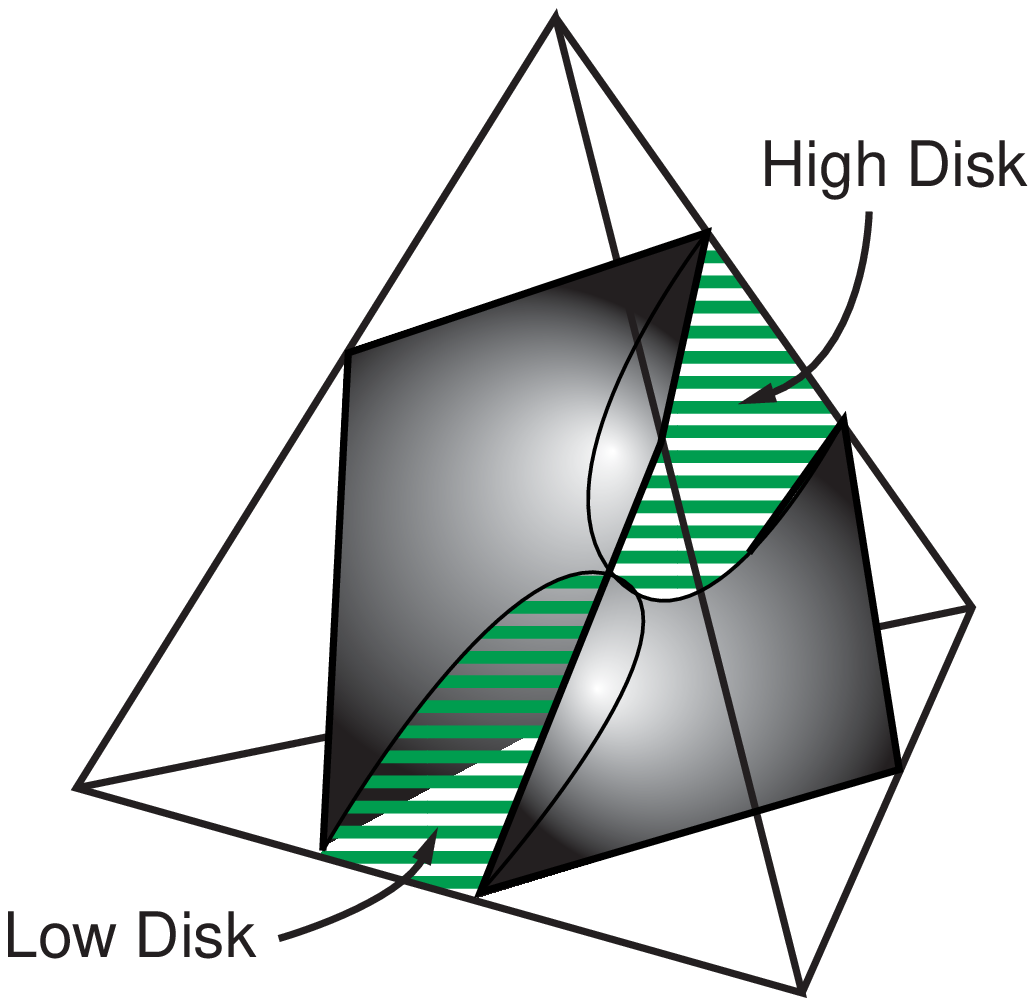}
        \caption{Relative compressing disks for an octagonal piece of $G_j$.}
        \label{f:Octagon}
        \end{center}
        \end{figure}

We conclude that $G_j$ is made up of all normal disks, with the exception of either a single disk with a boundary curve of length 8, OR a single place where there are two normal disks tubed together by an unknotted tube. This is the precise definition of an almost normal surface.

Our proof is complete by noting that there must be an octagonal disk or a tube {\it somewhere}, because $G_j$ is a thick level, and hence there is at least one compressing or relative compressing disk on both sides. If there were no tubes or octagons, then we would not have this.
\end{proof}

The next theorem is an application of the results of this section. This result is an important step in the proof of Rubinstein and Thompson's recognition algorithm for $S^3$.

\begin{thm}
\label{s3}
{\rm [Rubinstein \cite{rubinstein:93}, Thompson \cite{thompson:94}]} Let $T$ be a triangulation of some 3-manifold, $M$. If $B$ is a ball in $M$ such that $\partial B$ is the unique normal 2-sphere in $B$, and such that $B$ does not contain any of the vertices of $T^0$, then $B$ contains an almost normal 2-sphere with an octagon as its exceptional piece. 
\end{thm}

\begin{proof}
Inside $B$ the 1-skeleton is simply a 1-manifold. Choose an HST splitting, $\{G_j\}$, of $(B,T^1)$ which is minimal. We first claim that $\{G_j\}$ is strongly irreducible. If not, then we can perform one of the four operations described in the proof of Lemma \ref{l:untang} to produce an HST splitting which is smaller than $\{G_j\}$. 

We now call upon Theorem \ref{t:an} to isotope $\{G_j\}$ so that all thick levels are almost normal. Suppose $G_1$ is such a thick level. Let $S$ be the sphere component of $G_1$ which contains the exceptional piece, which is either an octagon or a tube. In the latter case $G_2$ is obtained from $G_1$ by compressing the tube, yielding two normal 2-spheres. As the only such normal 2-spheres are normally parallel to $\partial B$, the sequence $\{\partial B, G_3, G_4, ...\}$ is also an HST splitting of $(M,\Gamma)$, contradicting the minimality of $\{G_j\}$. We conclude then that the exceptional piece of $S$ was an octagon.
\end{proof}

\section{Finding almost normal surfaces.}

In this section we tie together the main results of this paper to show how one can prove that various topologically interesting surfaces can be made almost normal. 

\begin{dfn}
If $M$ is an irreducible submanfold of a triangulated 3-manifold then $S \subset M$ is a {\it maximal normal 2-sphere} in $M$ if there is no other normal 2-sphere which bounds a ball in $M$ which contains $S$. 
\end{dfn}

If $M$ is not homeomorphic to $B^3$ or $S^3$, and has a unique maximal normal 2-sphere, $S$, then we recall two basic facts about $S$ from \cite{rubinstein:93} or \cite{thompson:94}. First, $S$ bounds a ball, $B$, in $M$ which contains all of the vertices of $T^0 \cap M$. Second, there are no almost normal 2-spheres in $M$ that are disjoint from $B$.

\begin{thm}
\label{t:main}
Suppose $M$ ($\neq B^3$ or $S^3$) is an irreducible submanifold of a triangulated 3-manifold which is bounded by a normal or almost normal surface, and contains at most one maximal normal 2-sphere. If $\{F_i\}$ is a strongly irreducible HST splitting of $M$ then each thin level is isotopic to a normal surface, and each thick level is isotopic to an almost normal surface. 
\end{thm}

\begin{proof}
Let $S$ denote a maximal normal 2-sphere in $M$. $S$ bounds a ball, $B$ which contains all of the vertices of $T^0 \cap M$. Collapsing $B$ to a point, $v$, then turns $T^1$ into a properly embedded graph, $\Gamma$, whose unique vertex is $v$. Note that any surface which is disjoint from $v$ after the collapse can be identified with a surface which is disjoint from $B$ before the collapse.

Let $\{F_i\}$ be a strongly irreducible HST splitting of $M$. Lemma \ref{l:untang} implies that there is a strongly irreducible HST splitting, $\{G_j\}$, of $(M,\Gamma)$ such that $[\{G_j\}]=\{F_i\}$. By Lemma \ref{l:underlyingodd}, for each thick (thin) level of $\{F_i\}$ there is a thick (thin) level of $\{G_j\}$ whose non-$S^2$ components are parallel to it. Note that as $\{G_j\}$ is a HST splitting of $(M,\Gamma)$, which, by definition, is a HST splitting of $(M-Nbhd(v),\Gamma^1)$, we can also think of $\{G_j\}$ as a strongly irreducible HST splitting of $(M-int(B),T^1)$. 

By Theorem \ref{t:an} we may isotope $\{G_j\}$ so that each thin level is normal and each thick level is almost normal. Note that for each almost normal thick level, the exceptional piece is not contained in any 2-sphere component, because there are no almost normal 2-spheres in $M-int(B)$. Hence, every thick level of $\{F_i\}$ is isotopic to an almost normal surface.
\end{proof}

\begin{thm}
\label{hs}
{\rm [Rubinstein \cite{rubinstein:93}, Stocking \cite{stocking:96}]} In any triangulation of a closed, irreducible 3-manifold other than $S^3$, any strongly irreducible Heegaard surface is isotopic to an almost normal surface.
\end{thm}

\begin{proof}
Let $F$ be a strongly irreducible Heegaard surface for a closed, irreducible 3-manifold, $M$. Then $\{\emptyset, F, \emptyset\}$ is a strongly irreducible HST splitting of $M$ with $F$ as a thick level. By the results of \cite{rubinstein:93} or \cite{thompson:94}, $M$ contains a unique maximal normal 2-sphere. Hence, Theorem \ref{t:main} implies that $F$ can be made almost normal. 
\end{proof}

\begin{thm}
\label{sch}
{\rm [Schleimer \cite{schleimer:01}]} Suppose $M$ is a submanifold of a 3-manifold with triangulation, $T$, such that $M$ is bounded by two topologically, but not combinatorially, parallel normal surfaces. If, in addition, $M$ does not contain a normal 2-sphere then $M$ contains a boundary parallel almost normal surface.
\end{thm}

\begin{proof}
By assumption, $M$ is homeomorphic to $F \times I$, for some subsurface, $F$, of $\partial M$. The sequence, $\{F \times \{0\},F \times \{1/2\},F \times \{1\}\}$, is then a strongly irreducible HST splitting of $M$, with $F \times \{1/2\}$ as a thick level. Theorem \ref{t:main} now implies that $F \times \{1/2\}$ can be made almost normal.
\end{proof}

\begin{thm}
\label{nonsep}
If a triangulated irreducible 3-manifold contains a non-separating surface, then it contains a non-separating almost normal surface.
\end{thm}

\begin{proof}
Let $S$ denote a maximal normal 2-sphere in an irreducible 3-manifold, $N$. As in the proof of Theorem \ref{hs}, $S$ bounds a ball which we collapse to a point, transforming $T^1$ to a 1-vertex graph, $\Gamma$. 

It is well known that if a 3-manifold contains a non-separating surface then it contains an essential one. We may also assume that we have such a surface, $F$, that is disjoint from $S$, and that this surface can be normalized in $N \backslash S$. Let $F_-$ and $F_+$ denote two parallel copies of $F$. Finally, let $M$ denote the submanifold of $N$ cobounded by $F_-$ and $F_+$ which contains $S$. 

As $M$ is cobounded by surfaces which are essential in $M$ there is a strongly irreducible HST splitting, $\{F_i\}_{i=0}^n$, of $M$ such that $F_0=F_-$ and $F_n=F_+$. Suppose $F_p$ is a thick level of this splitting. Then $F_p$ separates $F_-$ from $F_+$ in $M$, but is non-separating in $N$. Theorem \ref{t:main} may now be applied to $M$ to show that $F_p$ can be made almost normal.
\end{proof}

\begin{dfn}
Suppose $G$ is a normal surface in a triangulated 3-manifold, and $F$ is any surface disjoint from $G$. Then we say $F$ {\it normalizes monotonically to} $G$ if there is a $K$-compression body, $(W;K)$, such that $W \cap T^1=K$, $F=\partial _+ W$, and $G=\partial _-W$. 
\end{dfn}

The following theorem characterizes when a normal surface compresses in terms of the existence of certain almost normal surfaces. In \cite{jb:02} we use this to give a new algorithm to recognize when a given normal surface is compressible.

\begin{thm}
\label{compressible}
Suppose $M$ is a submanifold of a triangulated irreducible 3-manifold which is bounded by a connected normal surface, and contains at most one maximal normal 2-sphere. If $\partial M$ is compressible in $M$ then there is an almost normal surface in $M$ which is topologically parallel to $\partial M$, and normalizes monotonically to $\partial M$.
\end{thm}

\begin{proof}
Let $\mathcal C$ denote a maximal collection of compressing disks for $\partial M$ in $M$. Let $P$ denote the surface obtained from $\partial M$ by compressing along all disks in $\mathcal C$. Let $N_1$ denote side of $P$ opposite $\partial M$ in $M$. Alter $P$ by throwing away any 2-sphere component which bounds a ball in $N_1$, and continue to denote this new surface as $P$. If $P=\emptyset$, then $M$ is a handlebody, and $\{F_0, F_1\}=\{\emptyset,\partial M \}$ is a strongly irreducible HST splitting of $M$. Otherwise $P$ is essential in $N_1$, and we can find a strongly irreducible HST splitting, $\{F_i\}_{i=0}^{2n}$, of $N_1$, where $F_0=\emptyset$ and $F_{2n}=P$. In this case, let $F_{2n+1}=\partial M$, so that $\{F_i\}_{i=0}^{2n+1}$ is a strongly irreducible HST splitting of $M$. In any case, $F_{2n+1}=\partial M$ is a thick level. Hence, we may apply Theorem \ref{t:main} to show that it is isotopic to an almost normal surface in $M$. 
\end{proof}

As a Corollary to Theorem \ref{compressible} we obtain:

\begin{cor}
\label{normalheegaard}
If $H$ is a Heegaard surface for a 3-manifold which is isotopic to a normal surface in the complement of a maximal normal 2-sphere then $H$ is isotopic to two, topologically parallel, almost normal surfaces.
\end{cor}

\begin{proof}
Let $H_0$ denote a normal representative of $H$. Then $H_0$ bounds two compression bodies, $W$ and $W'$. Both of these satisfy the hypotheses of Theorem \ref{compressible}, and hence, both contain almost normal surfaces which are topologically parallel to $H_0$. 
\end{proof}

The following result is important because it gives an exact characterization of handlebodies in terms of the existence of normal and almost normal surfaces. In \cite{jb:02} we use this to give a new algorithm to recognize when a given almost normal surface is a Heegaard surface. 

\begin{thm}
\label{t:handlebody}
Suppose $W$ is a submanifold of a triangulated 3-manifold such that $\partial W$ is almost normal, and such that $W$ does not contain a normal 2-sphere. Then $W$ is a handlebody if and only if there is a sequence of pairwise disjoint surfaces, $\{G_i\}_{j=0}^{2n+1}$, such that 
\begin{enumerate}
    \item $G_{2n+1}=\partial W$,
    \item $G_0=\emptyset$,
    \item for each even $j$, $G_j$ is normal,
    \item for each odd $j$, $G_j$ is almost normal, and normalizes monotonically to $G_{j-1}$ and $G_{j+1}$, and
    \item $G_i \le G_j$ for all $i<j$ (with respect to the ordering given in Definition \ref{d:order}).
\end{enumerate}
\end{thm}

\begin{proof}
The proof that the existence of such a sequence implies that $W$ is a handlebody is left as an exercise to the reader. The proof of the converse is similar to that of Theorem \ref{t:main}. Suppose that $W$ is a handlebody that satisfies the hypotheses. Then $\{\emptyset, \partial W\}$ is a strongly irreducible HST splitting of $W$. Lemma \ref{l:untang} implies that that there is a strongly irreducible HST splitting, $\{G_j\}_{j=0}^m$, of $(W,T^1)$ such that $[\{G_j\}]=\{\emptyset, \partial W\}$. This immediately implies that $G_0=\emptyset$, $G_m=\partial W$, and that for all $i<j$, $G_i \le G_j$. By Theorem \ref{t:an} we may isotope $\{G_j\}$ so that each surface of even index is normal, and each surface of odd index is almost normal. Furthermore, by definition there is a $T^1$-compression body between consecutive thick and thin level, so each almost normal thick level normalizes monotonically to the next (and previous) normal thin level. As $\partial W$ is assumed to be almost normal, we conclude that $m=2n+1$, for some $n$. 
\end{proof}

\bibliographystyle{alpha}

\end{document}